\documentclass[11pt,
]{amsart}
\usepackage{bbm,times,url,graphicx,anysize,amssymb,amsmath,amsthm,multicol}
\title{Rooted induced trees in triangle-free graphs}
\author{Florian Pfender}
\address{Universit{\"a}t Rostock\\
Institut f{\"u}r Mathematik\\
D-18055 Rostock, Germany}
\email{Florian.Pfender@uni-rostock.de}
\date{}

\newtheorem{theorem}{Theorem}

\newtheorem{corollary}{Corollary}

\begin{document}
\maketitle

\begin{abstract}
For a graph $G$, let $t(G)$ denote the maximum number of vertices in an induced subgraph of $G$ that is a tree. Further, for a vertex
$v\in V(G)$, let $t(G,v)$ denote the maximum number of vertices in an induced subgraph of $G$ that is a tree, with the 
extra condition that the tree must contain $v$.
The minimum of $t(G)$ ($t(G,v)$, respectively) over all connected triangle-free graphs $G$ (and vertices $v\in V(G)$) on $n$ vertices is 
denoted by $t_3(n)$ ($t_3^*(n)$).
Clearly, $t(G,v)\le t(G)$ for all $v\in V(G)$. In this note, we solve the extremal problem of maximizing $|G|$ for given
$t(G,v)$, given that $G$ is connected and triangle-free. We show that $|G|\le 1+\frac{(t(G,v)-1)t(G,v)}{2}$ and determine the unique extremal graphs.
Thus, we get as corollary that $t_3(n)\ge t_3^*(n)=\lceil \frac{1}{2}(1+\sqrt{8n-7})\rceil$, improving a recent result by
Fox, Loh and Sudakov.
\end{abstract}


All graphs in this note are simple and finite. For notation not defined here we refer the reader to Diestel's book~\cite{Di}.

For a graph $G$, let $t(G)$ denote the maximum number of vertices in an induced subgraph of $G$ that is a tree. The problem
of bounding $t(G)$ was first studied by Erd\H{o}s, Saks and S\'{o}s~\cite{ESS} for certain classes of graphs, one of them
being triangle-free graphs. Let $t_3(n)$ be the minimum of $t(G)$ over all connected triangle-free graphs $G$ on $n$ vertices.
Erd\H{o}s, Saks and S\'{o}s showed that
$$
\Omega\left(\frac{\log{n}}{\log\log{n}}\right)\le t_3(n)\le O(\sqrt{n}\log{n}).
$$
This was recently improved by Matou\u{s}ek and \u{S}\'{a}mal~\cite{MS} to
$$
e^{c\sqrt{\log{n}}}\le t_3(n)\le 2\sqrt{n}+1,
$$
for some constant $c$. For the upper bound, they construct graphs as follows.
For $k\ge 1$, let $B_k$ be the bipartite graph obtained from the path $P^k=v_1\ldots v_{k}$ if we replace $v_i$ by 
$\frac{k+1}{2}-|\frac{k+1}{2}-i|$ independent 
vertices for $1\le i\le k$. This graph has $|B_k|=\left\lfloor\frac{(k+1)^2}{4}\right\rfloor$ vertices,
yielding the bound.

For a vertex
$v\in V(G)$, let $t(G,v)$ denote the maximum number of vertices in an induced subgraph of $G$ that is a tree, with the 
extra condition that the tree must contain $v$. Similarly as above, we define $t_3^*(n)$
as the minimum of $t(G,v)$ over all connected graphs $G$ with $|G|=n$ and vertices $v\in V(G)$.
As $t(G,v)\le t(G)$ for every graph,
this can be used to bound $t_3(n)$. In a very recent paper, Fox, Loh and Sudakov do exactly that to show that
$$
\sqrt{n}\le t_3^*(n)\le t_3(n) ~\mbox{ and }~ t_3^*(n)\le \lceil \tfrac{1}{2}(1+\sqrt{8n-7})\rceil.
$$
For the upper bound, they construct graphs similarly as above.
For $k\ge 1$, let $G_k$ be the bipartite graph obtained from the path $P^k=v_0v_1\ldots v_{k-1}$ if we replace $v_i$ by $k-i$ independent 
vertices for $1\le i\le k-1$. This graph has $|G_k|=1+\frac{(k-1)k}{2}$ vertices, yielding the bound.

In this note, we show that this upper bound is tight, and that the graphs $G_k$ are, in a way, the unique extremal graphs. 
This improves the best lower bound on $t_3(n)$ by a factor of roughly $\sqrt{2}$. In~\cite{FLS}, the authors relax the problem
to a continuous setting to achieve their lower bound on $t_3^*(n)$. While most of our ideas are inspired by this proof,
we will skip this initial step and get a much shorter and purely combinatorial proof of our tight result.

\begin{theorem}\label{Tmain}
Let $G$ be a connected triangle-free graph on $n$ vertices, and let $v\in V(G)$. If $G$ contains no tree through $v$ on $k+1$ 
vertices as an induced subgraph, then $n\le 1+\frac{(k-1)k}{2}$. Further, equality holds only if $G$ is isomorphic to $G_k$ with $v=v_0$.
\end{theorem}
In the proof we will use the following related statement.
\begin{theorem}\label{Tside}
Let $G$ be a connected triangle-free graph, and let $v\in V(G)$. If $G$ contains no tree through $v$ on $k+1$ 
vertices as an induced subgraph, then $|V(G)\setminus N[v]|\le \frac{(k-2)(k-1)}{2}$. 
\end{theorem}
\begin{proof}[Proof of Theorems~\ref{Tmain} and~\ref{Tside}.]
Let $A(k)$ be the statement that Theorem~\ref{Tmain} is true for the fixed value $k$, and let $B(k)$
be the statement that Theorem~\ref{Tside} is true for $k$. We will use induction on $k$ to show $A(k)$ and $B(k)$ simultaneously. 

To start, note that $A(k)$ and $B(k)$ are trivially true for $k\le 2$. 
Now assume that $A(\ell)$ and $B(\ell)$ hold for all $\ell<k$ for some $k\ge 3$, and we will show $B(k)$.
We may assume that every vertex in $N(v)$ is a cut vertex in $G$ (otherwise delete it and proceed with the smaller graph). 
Let $N(v)=\{ x_1,x_2,\ldots,x_r\}$, and let $X_i$ be a component of $G\setminus N[v]$
adjacent  only to $x_i$ for $1\le i\le r$. 

Let $k_i+1$ be the size of a largest induced tree in $x_i\cup X_i$ containing $x_i$. 
Clearly, $G$ contains an induced tree through $v$ on $1+r+\sum k_i$ vertices, so $1+r+\sum k_i\le k$ (and in particular $k_i+1< k$).
By $A(k_i+1)$ we have $|X_i|\le \frac{k_i(k_i+1)}{2}$. 

Now replace each
$G[x_i\cup X_i]$ by a graph isomorphic to $G_{k_i}$ with $v_0=x_i$, reducing the total number of vertices by at most
$\sum k_i$. Note that this new graph $G'$ is triangle-free and connected. Since every maximal induced tree in $G$ through $v$ must contain a vertex $x_i$ for some $1\le i\le r$, and therefore exactly $k_i$ vertices of $X_i$, every induced tree through $v$ in $G'$ has fewer than $k$
vertices. Therefore, by $B(k-1)$, 
$$
|V(G)\setminus N[v]|\le |V(G')\setminus N[v]|+\sum k_i \le \frac{(k-3)(k-2)}{2}+k-r-1\le\frac{(k-2)(k-1)}{2},
$$
establishing $B(k)$.
Equality can hold only for $r=1$, and if $G[x_1\cup X_1]$ is isomorphic to $G_{k-1}$ by $A(k-1)$. Further,
every vertex in $N(v)$ must be adjacent to all neighbors of $x_1$ as otherwise a tree on $k+1$ vertices could be found in $G$. To see $A(k)$, note that $|N(v)|\le k-1$ or there is an induced star centered at $v$.
\end{proof}

As a corollary we get the exact value for $t_3^*(n)$, which is an improved lower bound for $t_3(n)$.
\begin{corollary}
$\lceil \frac{1}{2}(1+\sqrt{8n-7})\rceil= t_3^*(n)\le t_3(n)\le 2\sqrt{n}+1$.
\end{corollary}

\section*{Concluding remarks}
One may speculate that, similarly to the role of the $G_k$ for $t_3^*(n)$, the graphs $B_k$ are extremal graphs for $t_3(n)$. This is not true for $k=5$, though, as $K_{5,5}$ minus a perfect matching has no induced tree with more than $5$ vertices, and $B_5$ has only $9$ vertices. We currently know of no other examples beating the bound from $B_k$. In fact, with a similar proof as above one can show that $B_k$ is extremal under the added condition that $G$ has diameter $k-1$.

\bibliographystyle{amsplain}

\end{document}